\documentclass[english,a4paper,12pt]{article}
\usepackage{babel}
\usepackage[ansinew]{inputenc}
\usepackage{graphicx}
\usepackage{verbatim}
\usepackage{amsmath}
\usepackage{amssymb}
\def\C{\mathbb C}
\def\Cd{\widehat{\mathbb C}}
\def\N{\mathbb N}

\def\Z{\mathbb Z}

\newtheorem{theo}{Theorem}[section]
\newtheorem{cor}[theo]{Corollary}
\newtheorem{ex}[theo]{Example}
\newtheorem{lemma}[theo]{Lemma}
\newtheorem{prop}[theo]{Proposition}
\newtheorem{defi}[theo]{Definition}
\newtheorem{rem}[theo]{Remark}

\def\pr{{\it Proof. }}
\begin{document}

\title{Rational functions that share finite values with their first derivative}
\author{Andreas Sauer* \and Andreas Schweizer**}
\date{\empty}
\maketitle
\thispagestyle{empty}
\textit{\small * Hochschule Ruhr West, Duisburger Str.~100, 45479 M\"{u}lheim an der Ruhr,\\ Germany, andreas.sauer@hs-ruhrwest.de}\\[.3cm]
\textit{\small ** Department of Mathematics Education, Kongju National University, Gongju 32588,\\South Korea, schweizer@kongju.ac.kr}

\begin{abstract}
We treat shared value problems for rational functions $R(z)$ and their derivative $R'(z)$ in the plane and on the sphere. We also consider shared values for the pair $R(w)$ and $\partial_{z} R = \lambda w \cdot R'(w)$ on $\C \setminus \{ 0\}$ and $\Cd$, again with rational functions $R$. In $\C \setminus \{ 0\}$ this is related to shared values of meromorphic functions $f \colon \C \to \Cd$ and $f'$ through $f(z) = R(w)$ with $w = \exp(\lambda z)$, while on $\Cd$ this is connected to shared limit values in a similar fashion.
\end{abstract}
\renewcommand{\thefootnote}{}
\footnotetext{\hspace*{-.51cm}2020 Mathematics Subject Classification: Primary: 30D35, Secondary: 30C10\\ %
}
\tableofcontents
\noindent
\section{Introduction}
In this paper we will investigate shared value problems for rational functions and their first derivative in various relevant settings.\\

To the best of our knowledge, this has not been treated systematically in the past. On the other hand, it is close to impossible to check the considerable body of work on shared value problems, in order to be absolutely sure about the originality of all the statements and examples in this paper. If we missed appropriate references, we would like to apologize in advance.\\

There are several motivations for this paper, and we will mention these in the appropriate sections. Nonetheless, we want to start here with the main cause for our investigations:\\

Let $f \colon \C \to \Cd$ be a meromorphic function, including transcendental functions. A value $a \in \C$ is a {\it shared} value of $f$ and $f'$ if $f(z) = a \Leftrightarrow f'(z) = a$ for all $z \in \C$. Note that $f$ and $f'$ always have the same poles, so that $\infty$ is shared automatically. If, in addition, in every $a$-point the order of $f$ is equal to the order of $f'$, then $a$ is said to be shared {\it counting multiplicities} (CM). For $a \in \C^* := \C \setminus \{0\}$ the assumption CM results in $f$ and $f'$ having only simple $a$-points, while for $a=0$ this means, that $f$ and $f'$ have no zeros at all.\\

Shared values for $f$ and $f'$ is a well-known topic in uniqueness theory for meromorphic functions, we refer to \cite{YaYi}, section 8. The most famous results in this respect are the following theorems.
\\ \\
{\bf Theorem A.} \cite[Theorem 3]{Gu}, \cite[Satz 2]{MuSt1}  \it
Let $f \colon \C \to \Cd$ be a non-constant meromorphic function that shares three finite 
values with $f'$, then $f\equiv f'$.
\rm
\\
\\
{\bf Theorem B.} \cite[SATZ, p.~49]{MuSt2}  \it
Let $f \colon \C \to \Cd$ be a non-constant meromorphic function that shares two finite values CM with $f'$, then $f\equiv f'$.
\rm
\\ \\
Both theorems can be found, including proofs, in \cite{YaYi}, Theorem 8.22 and 8.20, respectively.\\

It is known that the assumption CM in Theorem B cannot be dropped completely, as there are two examples where $f$ and $f'$ share two finite values, but $f \not \equiv f'$ (see Example \ref{example_R_dR}). Both of these examples have in common, that $0$ is one of the shared values.\\ 

It is a long standing open problem, whether the assumption CM in Theorem B can be weakened, if the two shared values are both non-zero (see \cite{MuSt2}, BEISPIEL B and the comments thereafter). As far as we are aware of, the only progress that has been made concerning this question are the partial results in \cite{Zh}, \cite{Li} and \cite{Ta}. So far, it is not even known, whether there exists a rational function $f$, such that $f$ and $f'$ share two (non-zero) values. Note that such an example would be a counter-example to $f \equiv f'$ since this is impossible for rational functions other than $f \equiv 0$.\\

In this paper we prove in this connection:\\

\noindent
{\bf Theorem \ref{rational_plane_IM_ab}} \it A non-constant rational function $f \colon \C \to \Cd$ cannot share two finite values with its derivative. 
\rm
\\ \\
Note that in this theorem one of the shared values is allowed to be zero.\\

Like most relevant examples in uniqueness theory, the above mentioned two examples in Example \ref{example_R_dR} are of the form $f(z) = R \left( \text{e}^{\lambda z} \right)$ with a rational function $R$ and $\lambda \neq 0$ (with $\deg(R) \le 2$). Again in connection with the problem from \cite{MuSt2} we prove:\\

\noindent
{\bf Theorem \ref{R_dR_aCMbIM_punctured_plane} \it Let $R$ be a non-constant rational function, $\lambda \neq 0$ and $f(z) = R \left( \text{\rm e}^{\lambda z} \right)$. If $f$ and $f'$ share two finite values $a \neq 0$ and $b$, such that $a$ is shared CM, then $f \equiv f'$.
\rm
\\ \\
We describe later, what $f \equiv f'$ means in terms of $R$.\\

We also treat the case where $a=0$ is shared CM by $f$ and $f'$ with this form, which is another way of saying that $f$ omits the value $0$ (see Theorem \ref{R_dR_a=0CMbIM_punctured_plane}), but of course we get the same conclusions as in \cite{Zh}.\\

As a last remark on this topic we want to mention that we prove in Proposition \ref{prop_deg2_punctured_plane}, that for $\deg(R) \le 2$ no other such examples exist than the ones given in Example \ref{example_R_dR}. Any new non-trivial example of this form would have $\deg(R) \ge 3$.\\

\noindent
Now we will focus on rational functions, which we mostly denote by $R$, and we fix some notation. We consider rational functions $R$ as holomorphic mappings from the Riemann sphere $\Cd$ into $\Cd$. We generally exclude the function $R \equiv \infty$. Let $R=P/Q$ be a representation of a non-constant rational function $R$ by the quotient of two relatively prime polynomials $P$ and $Q$. It is an easy consequence of the fundamental theorem of algebra that the preimages in $\Cd$ for any $a \in \Cd$, when counted with regard to multiplicities, give the same constant $d$, the {\it degree} of $R$, independently of $a$. The degree $\deg(R)$ of $R$ will always be denoted by $d$.\\

$R$ and $R'$ {\it share} a value $a \in \Cd$ in a set $G \subset \Cd$ if the equations $R(z) = a$ and $R'(z) = a$ have the same solutions $z \in G$, where $z = \infty$ is possible if $\infty \in G$. If, in addition, in every $a$-point the order of $R$ is equal to the order of $R'$, then $a$ is said to be shared {\it counting multiplicities} (CM), where the order of $R$ at $z=\infty$ is defined in the usual way for points on Riemann surfaces. If we want to emphasize that we do not assume CM, then we use the standard notation IM for {\it ignoring multiplicities}. Note that if $R$ has a pole of first order at $z=\infty$, then $R'$ has a non-zero finite value at $z=\infty$, and if $R$ has a finite value at $z=\infty$, then $R'(\infty)=0$. The latter statement gives no information on the order of $R$ at $z=\infty$, though.\\

Let $r$ be the total ramification of $R \colon \Cd \to \Cd$, i.e.~the sum $\sum_{z \in \Cd} (\text{mult}_f(z)-1)$, where $\text{mult}_f(z)$ denotes the multiplicity (or order) of $f$ at the point $z \in \Cd$. On several occasions we will use the Riemann-Hurwitz formula, which says that
$$
r = 2d-2.
$$
The Riemann-Hurwitz formula is often proved for arbitrary compact Riemann surfaces, but we do not need this general form. Again let $R=P/Q$ be a representation by the quotient of two relatively prime polynomials $P$ and $Q$. With the help of the fundamental theorem of algebra it is an easy exercise to prove $r = 2d-2$ by counting according to multiplicities the zeros of $R' = (P'Q-PQ')/Q^2$ and the poles of $R$, taking into account the different possibilities of how $R$ can behave at $\infty$.\\

Another argument that we will use, is a version of the residue theorem. Again, there exist statements for compact Riemann surfaces. But we just need the elementary, yet powerful, formulation for rational functions: The sum of the residues at poles of $R$ in $\C$ equals the negative of the coefficient of $1/z$ in the Laurent expansion of $R$ at $z = \infty$. This follows from the classical residue theorem for meromorphic functions in $\C$ and a simple application of the transformation rule to line integrals, taking into account that $R(1/z)$ is a rational function.\\

The necessary basic notions and techniques from complex analysis to prove the aforementioned facts on rational functions can be found for instance in \cite{A}.

\section{Polynomials}
The case of polynomials $P$ (with degree $d$) is particularly simple. It is a matter of taste, whether one considers $P$ as a function on $\C$ or $\Cd$. If $d$ is at least $2$ then $P$ and $P'$ have a single common pole at $z = \infty$. Note that values cannot be shared CM since $\deg(P') = d - 1$.
\begin{theo} \label{polynomial_non_zero} A non-constant polynomial cannot share a finite non-zero value with its derivative.
\end{theo}

\pr Assume that $P$ and $P'$ share $a \in \C^*$ and consider $\Psi := (P'-a)/(P-a)$. Then $\Psi$ is a rational function without any poles in $\Cd$. Therefore $\Psi$ is constant. Since $\Psi(\infty)=0$ we conclude $\Psi \equiv 0$ which implies $P' \equiv a$. Since $a$ is shared by $P$ and $P'$ this gives the contradiction $P \equiv a$. $\square$ \\

If $P$ and $P'$ share two finite values at least one of them is non-zero. This immediately gives the following consequence.
\begin{cor} \label{polynomial_two_finite} A non-constant polynomial cannot share two finite values with its derivative.
\end{cor}

Hence we are left with the case that $P$ and $P'$ share the value $a=0$.
\begin{theo} \label{polynomial_zero} A non-constant polynomial $P$ shares the value $0$ with its derivative if and only if $P(z) = c \cdot (z-A)^n$ with arbitrary $A$, $c \neq 0$ and $n \ge 2$.
\end{theo}
\pr  Let $n_0$ be the number of zeros of $P$ in $\C$, counted without regard to multiplicities. Since $P'$ is allowed to have zeros only at points where $P$ is zero, and since each such zero of $P'$ has an order that is one less than the order of $P$, we get the equation $\deg(P') = d - n_0$. Since for polynomials always $\deg(P') = d - 1$ it follows $n_0 = 1$. The rest is easy. $\square$
\section{Rational functions on the sphere} \label{Rational_functions_sphere}
Now we consider rational functions $R \colon \Cd \to \Cd$ and the values have to be shared on the whole sphere. Shared values on $\Cd$ were already considered in \cite{Pi}, and this section can be seen as a supplement to \cite{Pi}, as well as for \cite{Sa_compact} and \cite{Sch}.

\begin{theo} \label{rational_zero_implied} Let $R \colon \Cd \to \Cd$ be a non-constant rational function. Then $R'(z) = 0 \Rightarrow R(z)=0$ for all $z \in \Cd$ if and only if $R$ has the form $R(z) = c \cdot (z-A)^n$ with arbitrary $A$, $c \neq 0$ and $n \in \Z \setminus \{ 0\}$.
\end{theo}

\pr \glqq$\Rightarrow$\grqq. From the assumption it follows that no finite non-zero value is ramified for $R$. Note that $R(\infty)$ cannot be finite and non-zero since then $R'(\infty)=0$ in contradiction to the assumption. Hence the total ramification of $R$, which is $2 d - 2$, is contained in the zeros and the poles of $R$. It follows that $R$ has exactly one pole and one zero, both of order $d$. As already noted we have $R(\infty)=0$ or $R(\infty)=\infty$. This gives $R(z) = c \cdot (z-A)^n$ with arbitrary $A$, $c \neq 0$ and $n \in \Z \setminus \{ 0\}$.\\
\glqq$\Leftarrow$\grqq. This is easy to check. $\square$ \\

For $R(z) = c \cdot (z-A)^n$ with $n=1$ clearly $R'$ is a non-zero constant and for $n \in \Z \setminus \{ 0, 1\}$ these $R$ share $0$ with $R'$. Hence we immediately get the following corollary.
 
\begin{cor} \label{rational_zero_shared} A non-constant rational function $R \colon \Cd \to \Cd$ shares the value $0$ with its derivative if and only if $R(z) = c \cdot (z-A)^n$ with arbitrary $A$, $c \neq 0$ and $n \in \Z \setminus \{ 0, 1\}$.
\end{cor}

The next theorem is one of the main results of this section. It will also be a useful argument in later proofs. 

\begin{theo} \label{rational_two_values} A non-constant rational function $R \colon \Cd \to \Cd$ cannot share two finite values with its derivative.
\end{theo}

\pr First we assume that the shared values $a$ and $b$ are both non-zero. Note that $R(\infty) \neq a$ and $R(\infty) \neq b$, since otherwise $R'(\infty) = 0$, so that $a$ or $b$ would not be shared. It follows that all $a$ and $b$-points of $R$ and $R'$ lie in $\C$ and are simple for $R$. We consider
$$
\Psi := \frac{R' - R}{(R-a)(R-b)}.
$$
It is easy to see that $\Psi$ has no poles in $\Cd$. Therefore $\Psi$ is constant. $\Psi \equiv 0$ is impossible since then $R'=R$. Hence we get
$$
R' = R + c(R-a)(R-b)
$$
with $c \neq 0$. From this it follows that $\deg(R') = 2 d$, since the right hand side is a polynomial of degree two in $R$. Hence $\deg(R') > d$, so at least one of the $a$-points of $R'$ is ramified, because their number is $d$. This means that there exists $z_a \in \C$, such that $R'(z_a) = a$ and $R''(z_a)=0$. This and the rule of de L'Hospital show
$$
\Psi(z_a) = \frac{R''(z_a)-R'(z_a)}{R'(z_a)(R(z_a)-a) + R'(z_a)(R(z_a)-b)} = \frac{1}{b-a}.
$$
For exactly the same reason there exists $z_b \in \C$, such that $R(z_b) = b$ and $R''(z_b)=0$. This gives
$$
\Psi(z_b) = \frac{R''(z_b)-R'(z_b)}{R'(z_b)(R(z_b)-a) + R'(z_b)(R(z_b)-b)} = \frac{1}{a-b}.
$$
But this contradicts that $\Psi$ is constant, which proves the claim for non-zero $a$ and $b$.

Now assume that $a=0$ and $b \neq 0$ are the shared values. From Corollary \ref{rational_zero_shared} we know that $R$ has the form $R(z) = c \cdot (z-A)^n$ with $n \in \Z \setminus \{ 0, 1\}$. The $b$-points of $R$ are $1$-points of $R'/R$ and from the form of $R$ it follows $\deg(R'/R) = 1$. Thus $R$ has at most one $b$-point, which is not ramified since $R'= b \neq 0$ at this point, and hence $d = 1$. It follows that $R(z) = c \cdot (z-A)^{-1}$. It is easy to check that in this case $R'$ has two $b$-points, so that $b$ cannot be shared. $\square$ \\

We now treat the remaining problem, namely that $R$ and $R'$ share one non-zero finite value $a$. The following Lemma is a consequence of Theorem 1 in \cite{Gu}.

\begin{lemma} \label{gundersen_representation} Let $R \colon \Cd \to \Cd$ be a non-constant rational function that shares the finite non-zero value $a$ with its derivative $R'$ on the sphere. Then there exists a constant $c \in \C$ and a polynomial $P$, with zeros that are all simple, such that
$$
R = a \left( 1 + \frac{P}{cP + P'} \right) \qquad \text{and} \qquad R' = a \left( 1 - \frac{P[c(cP+P') + cP' + P'']}{(cP + P')^2} \right).
$$
The $a$-points of $R$ and $R'$ are exactly at the zeros of $P$. If $c(cP+P') + cP' + P''$ has a zero at $z_0$, then either $P(z_0)=0$ or $(cP + P')(z_0)=0$. If at such a point $z_0$ we have $P(z_0)=0$ then $z_0$ is an $a$-point of $R$ and $R'$ where the multiplicity of $R'$ is at least $2$ and therefore higher than the multiplicity of $R$ which is $1$. If at such a point $(cP + P')(z_0)=0$, then $z_0$ is a pole of $R$ which is at least of order $2$.
\end{lemma}

\pr Since $R$ and $R'$ share $a$ in the plane, Theorem 1 in \cite{Gu} shows that there exists an entire function $h$ with only simple zeros such that $R = a(1 + h/h')$. Because $R$ is rational, $h$ has finitely many zeros. Hence there exists a polynomial $P$ such that $h = P \exp(\alpha)$ with an entire function $\alpha$. Since 
\begin{equation} \label{alpha}
\frac{a}{R-a} = \frac{P'}{P} + \alpha'
\end{equation}
it follows that $\alpha$ is a polynomial. This gives $R = a(1 + P/(\alpha'P + P'))$. If the degree of $\alpha'$ is at least $1$, then $R(\infty)=a$ which implies $R'(\infty)=0$, so that $a$ would not be shared on the sphere. Hence $\alpha' \equiv c \in \C$. The form of $R'$ follows by simple calculations and all other statements are direct consequences of the representations of $R$ and $R'$. We only mention that if $c(cP+P') + cP' + P''$ has a zero at $z_0$ of order $n$ that is a zero of $cP + P'$ (of order $n+1$), then the order of the pole at $z_0$ of $R'$ is $2(n+1) - n = n + 2 \ge 3$. This implies that $R$ has a pole at $z_0$ at least of order $2$. $\square$

\begin{rem} \label{rem_cP}\rm In all examples of the form $R = a \left( 1 + P/(cP + P') \right)$, that we were able to find, we have $c=0$, i.e.~$R = a \left( 1 + P/P' \right)$. It seems to us that if $c(cP+P') + cP' + P''$ has zeros only at zeros of $P$ or $cP+P'$, then it follows $c=0$, but we were not able to find a proof for this statement. 
\end{rem}

\noindent
Unlike the case of polynomials it is possible that a general rational function shares a finite value CM with its derivative on the sphere. The examples can be determined completely and form a rather simple family of functions.

\begin{theo} \label{rational_CM} A non-constant rational function $R \colon \Cd \to \Cd$ shares the finite value $a$ CM with its derivative if and only if $a \neq 0$ and
$$
R(z) = a \left( 1 + \frac{(z-p)^{n+1} + C}{(n+1)(z-p)^{n}} \right)
$$
with $n \in \N$, $p \in \C$ and $C \in \C^*$.
\end{theo}

\pr Since $a$ is shared CM it follows $d = \deg(R')$. This is only possible if $R$ has exactly one pole $z=p$ in $\C$ and a pole at $z=\infty$. Thus $R$ is not of the form given in Corollary \ref{rational_zero_shared}, and therefore $a \neq 0$. Since $a$ is shared CM we conclude from Lemma \ref{gundersen_representation} that if $c(cP + P') + cP' + P''$ is zero at $z_0$ then $(cP + P')(z_0) = 0$. By Taylor expansion it is easy to see that the order of the zero $z_0$ of $cP + P'$ is exactly $+1$ higher than the order of $c(cP + P') + cP' + P''$. Hence the degree of $cP + P'$ is higher than the degree of $c(cP + P') + cP' + P''$, which is possible only if $c=0$. Since $R$ has only one pole in $\C$ it follows that $P'$ has only one zero in $\C$, so that $P'(z) = K_1 (z-p)^n$ with $K_1 \neq 0$, and therefore $P(z) = (K_1/(n+1)) (z-p)^{n+1} + K_2$ and $R$ is as stated in the theorem. Note that $C = (n+1)K_2/K_1 \neq 0$ is necessary, since otherwise $R$ has no pole in $\C$. Easy calculations show that such $R$ indeed share the value $a$ CM with their derivatives. $\square$

\begin{rem} \rm By simple transformations one can assume without loss of generality that $p=0$ and $a=1$. So in some sense
$$
R(z) = 1 + \frac{z^{n+1} + C}{(n+1)z^{n}}
$$
with $C \neq 0$ and $n \in \N$ are the only rational functions that share a finite value CM with their derivatives on $\Cd$. The shared value $1$ is taken by $R$ and 
$$
R'(z) = 1 - \frac{n \left( z^{n+1} + C \right)}{(n+1)z^{n+1}}
$$
in the $(n+1)$st roots of $-C$.
\end{rem}

A shared value $a$ for two functions $f$ and $g$ is said to be shared by {\it different multiplicities} (DM), if at each preimage of $a$ the order of $f$ and $g$ are different. If we assume that $R$ and $R'$ share $a \in \C^*$ DM, then all $a$-points of $R$ are simple, while $R'$ takes the value $a$ with multiplicity greater than one.

\begin{theo} \label{rational_DM} A non-constant rational function $R \colon \Cd \to \Cd$ cannot share a finite non-zero value $a$ by DM with its derivative.
\end{theo}

\pr Suppose $a$ is shared by DM by $R$ and $R'$. Since $\deg(R') \le 2d$ it follows immediately that $\deg(R') = 2d$ and that all $a$-points of $R'$ are double. With the above notation we conclude that all zeros of $c(cP + P') + cP' + P''$ are zeros of $P$ and also all zeros of $P$ are zeros of $c(cP + P') + cP' + P''$. Since all of these zeros are simple it follows that $c(cP + P') + cP' + P'' = c^2 P$ which implies $2cP' + P'' \equiv 0$ which gives the contradiction $P' \equiv 0$. $\square$ \\

It is possible though, that a non-zero value is shared by $R$ and $R'$ neither by CM, nor by DM. This is shown by the next example.

\begin{ex} \label{rational_IM} \rm The following examples $R$ share the value $a$ with $R'$ on the sphere, but not CM:
$$
R(z) = a \left( 1 + \frac{z^{n+1} + Cz}{(n+1)z^{n}+C} \right)
$$
with $C \neq 0$ and $n \ge 2$. (For $n=1$ we get one of the CM examples from Theorem \ref{rational_CM}.) The shared value $a$ is taken by $R$ and $R'$ in the $n$th roots of $-C$ and at $z=0$. The $a$-points of $R$ and $R'$ in the $n$th roots of $-C$ are simple for $R$ and $R'$, while the $a$-point of $R$ in $z=0$ is simple and for $R'$ it has order $n$. We were not able to find examples that are not of the form as given in Theorem \ref{rational_CM} or this example.
\end{ex}

\section{Rational functions in the plane} \label{Rational_functions_plane}
We now assume that $R$ and $R'$ share values with respect to $z \in \C$.\\

This setting is closely connected to classical shared value problems for meromorphic functions $f \colon \C \to \Cd$, we simply exclude that $f$ can be transcendental. Hence known theorems for general meromorphic functions $f$, which we have mentioned in the introduction, are true for $R$. For example, a non-constant rational function $R \colon \C \to \Cd$ cannot share three finite values with its derivative, which follows directly from Theorem A. The purpose of this section is to improve on these theorems for rational functions.\\

It is an interesting question whether a non-constant rational function $R$ exists, such that $R$ and $R'$ share two finite non-zero values in $\C$. If such an example existed, it would give an answer to the old question for meromorphic functions $f \colon \C \to \Cd$: {\it If $f$ and $f'$ share two finite non-zero values, does this imply $f \equiv f'$?} (see the introduction). We will prove that no rational counterexamples for this question exist. In fact, there are no such rational functions, even if one of the shared values is $0$.

\begin{theo} \label{rational_plane_IM_ab} A non-constant rational function $R \colon \C \to \Cd$ cannot share two finite values with its derivative.
\end{theo}

\pr Let $a$ and $b$ be the two finite shared values.\\

\noindent {\bf A)} First we assume that $a$ and $b$ are both non-zero.\\

If $R$ is a polynomial then the theorem follows from Theorem \ref{polynomial_two_finite}, and hence we can assume that $R$ has a pole in $\C$. If $R(\infty)$ and $R'(\infty)$ are different from $a$ and $b$, then $a$ and $b$ are shared on $\Cd$ and the claim follows from Theorem \ref{rational_two_values}. We have to consider two cases:\\

{\bf A i)} First we assume (without loss of generality) that $R(\infty)=b$ with order $m$. We prove $m=1$. Let $n_{\infty}$ be the number of different poles of $R$ in $\C$ and set
$$
\Psi := \frac{R' - a}{R - a}.
$$
Then $\Psi$ has $n_{\infty}$ simple poles with negative integer residues since $R'/(R - a)$ is a logarithmic derivative and $a/(R - a)$ is holomorphic at poles of $R$. It is elementary to check that $\Psi(\infty)= -a/(b-a)$ with order $m$. The residue theorem implies that the residue of $\Psi$ at $z = \infty$ equals the negative of the sum of the residues in $\C$, in particular it is non-zero. This is only possible if $m=1$.\\

It follows that $R$ has $d-1$ simple $b$-points in $\C$ which are $1$-points of $\Psi$. Since by counting the poles it holds $\deg(\Psi) = n_{\infty}$ we get $n_{\infty} \le d \le n_{\infty} + 1$. Furthermore, since $z = \infty$ is not a pole of $R$, we have $\deg(R') = d + n_{\infty}$. Hence only two cases are possible, either $\deg(R') = 2d$ or $\deg(R') = 2d -1$.\\

We now show $d \ge 4$ by estimating the ramification of $R'$. At $z=\infty$, which is a double zero of $R'$, the ramification is $1$, and in the poles it is $d + n_{\infty} - n_{\infty} = d$. In the $d-1$ $b$-points we have a ramification of $d + n_{\infty} - (d-1) = n_{\infty} +1$, and in the $d$ $a$-points it is $d + n_{\infty} - d = n_{\infty}$. Hence the ramification of $R'$ on $\Cd$ is at least $d + 2 n_{\infty} + 2$, while the total ramification of $R'$ is $2(d + n_{\infty}) - 2$, and therefore $d + 2 n_{\infty} + 2 \le 2 d + 2 n_{\infty}- 2$. This is equivalent to $d \ge 4$.\\

Next we consider
$$
\Phi :=  \frac{R' - R}{(R - a)(R-b)}.
$$

Since $a$ and $b$ are shared in $\C$ it follows that $\Phi$ has no poles in $\C$. The only pole of $\Phi$ at $z = \infty$ is of order $m$ since  $(R' - R)/(R - a)$ equals $-b/(b-a)$ at $z=\infty$. We have already proved $m=1$ and we conclude that $\Phi$ is a polynomial of degree one.\\

Since $d \ge 4$ clearly $\deg(R') \ge 2d-1 > d$. Hence there are $a$- and $b$-points in $\C$ where $R'$ is ramified. Just like in the proof of Theorem \ref{rational_two_values} we conclude that at a ramified $a$-point $z_a$ of $R'$ it holds $\Phi(z_a) = 1/(b-a)$ and at a ramified $b$-point $z_b$ of $R'$ it holds $\Phi(z_b) = 1/(a-b)$. Since $\Phi$ is linear (and non-constant) we conclude that there is exactly one such $z_a$ and exactly one such $z_b$. Considering $\Phi' / \Phi$ and using de L'Hospital's rule it is elementary to check that if the order of $R'$ at $z_a$ is greater or equal to $3$, then $\Phi'(z_a) = a/(a-b)^2$ and equally if the order of $R'$ at $z_b$ is greater or equal to $3$, then $\Phi'(z_b) = b/(a-b)^2$. If the degree of $R'$ is $2d$ then the order of $R'$ at $z_b$ is $d+2$ and at $z_a$ the order of $R'$ is $d+1$. Since $\Phi'$ is constant, it is not possible that the orders in $z_a$ and $z_b$ are both at least $3$. We get $d + 1 \le 3$, and hence $d \le 2$. If the degree of $R'$ is $2d-1$ then the order of $R'$ at $z_b$ is $d+1$ and at $z_a$ the order of $R'$ is $d$. Again since $\Phi'$ is constant we get $d \le 3$. Both possibilities contradict $d \ge 4$. This proves the theorem in case A i).\\ 

{\bf A ii)} In the second case we assume $R'(\infty) = b$. Then $R$ has a pole of order one in $z = \infty$ and it follows immediately $\Phi(\infty) = 0$. Since $\Phi$ is a polynomial we get $\Phi \equiv 0$, and hence $R' \equiv R$ which is impossible for non-constant rational $R$.\\

\noindent
{\bf B)} Now we assume $a = 0$ and $b \neq 0$.\\

If $R$ is a polynomial then the theorem follows from Theorem \ref{polynomial_two_finite}, and hence we can assume $\deg(R') \ge d$. If $R$ has no zero in $\C$ then $R(\infty) = R'(\infty) = 0$ and $0$ is shared on $\Cd$. By Corollary \ref{rational_zero_shared} we have $R(z) = c \cdot (z - A)^{-n}$ with $n \in \N$. It is easy to see, that such $R$ cannot share $b \neq 0$ with $R'$ in $\C$. Hence we can assume that $R$ has zeros in $\C$ and we denote their number by $n_0$. Since $0$ is shared it follows that $R'$ has $n_0$ less zeros than $R$ in $\C$ when counted according to multiplicities. From $\deg(R') \ge d$ we deduce that $R'(\infty)=0$. Therefore $R(\infty) \in \C$, where $R(\infty)=0$ is a contradiction to Corollary \ref{rational_zero_shared}. We set $K = R(\infty) \in \C^*$ with order $m$ and we denote by $n_{\infty} > 0$ the number of poles of $R$ in $\C$. Counting the zeros and the poles of $R'$ gives the equality $d - n_0 + m + 1 = d + n_{\infty}$ or $n_0 + n_{\infty} = m + 1$.\\

In order to show $K \neq b$ we assume $K=b$ and estimate the ramification $r$ of $R'$: In $z = \infty$ it is $m$. Since all poles lie in $\C$ and there are $n_{\infty}$ of them, we conclude that the ramification of $R'$ in the poles is $d + n_{\infty} - n_{\infty} = d$. In the $d - m$ $b$-points the ramification is $d + n_{\infty} - (d - m) = n_{\infty} + m$ and in the $n_0$ zeros it is $d + n_{\infty} - (m+1) - n_0$. Summing this up we conclude that $r$ is at least $r \ge m + d + n_{\infty} + m + d + n_{\infty} - (m+1)- n_0  = 2d + 2 n_{\infty} + m - 1 - n_0$. From $n_0 + n_{\infty} = m + 1$ we get $r \ge 2d + 3 n_{\infty} - 2$. Since the total ramification of $R'$ is $r = 2 (d + n_{\infty}) - 2 = 2d + 2 n_{\infty} - 2$ we obtain $n_{\infty}=0$. But then $R$ has no poles, and $R$ is constant.\\

Hence $K \neq b$ and $R$ has $d$ simple $b$-points in $\C$. We consider $L := R/R'$. Then $L$ has no poles in $\C$ and is a polynomial of degree $m+1$. $L$ has $d$ $1$-points in the $b$-points of $R$, which gives $d \le m+1$. As in A i) we apply the residue theorem to $\Psi := (R' - a)/(R - a)$ to obtain $m=1$. This implies $d=2$, $n_0=1$ and $n_{\infty}=1$, so that $R$ has exactly one zero and one pole in $\C$ both of order two. By translation we can assume that the pole of $R$ is in $z=0$, and by scaling $R$, $R'$ and $b$ we can assume $K=1$. Let $z=h$ be the double zero of $R$. Hence we can assume $R(z) = (z-h)^2/z^2$, so that $R'(z) = 2h(z - h)/z^3$ and $R''(z) = 2h(3h - 2z)/z^4$. It follows that $R''$ has exactly one zero in $\C$ at $z_0 = 3h/2$. Since $\deg(R')=3$ we conclude that $R'$ has at least one $b$-point that is ramified, and hence this $b$-point is in $z_0$. We get $b = 1/9 = R(z_0) = R'(z_0) = 8/(27h)$, and hence $h = 8/3$ and $z_0 = 4$. With these values the second $b$-point of $R$ is in $z_1=2$, but $R'(2) = -4/9 \neq b$. $\square$ \\

We now turn our attention to the case of one shared value. If the value is supposed to be shared by CM, then the conclusion of the next theorem is almost identical to Theorem \ref{rational_CM}, but now it is possible that the shared value is $a=0$. It is clear that in this case $a$ is omitted in $\C$ by $R$ and $R'$.

\begin{theo} \label{rational_CM_plane} A non-constant rational function $R \colon \C \to \Cd$ shares the finite value $a$ CM with its derivative if and only if
$$
R(z) = \frac{a}{n+1} (z-p) + a + \frac{C}{(z-p)^{n}}
$$
with $n \in \N$, $p \in \C$ and $C \in \C^*$.
\end{theo}

\pr \glqq$\Rightarrow$\grqq. We consider
$$
\Psi := \frac{R'-a}{R-a}.
$$
Since $a$ is shared CM, $\Psi$ has no zeros in $\C$ and has poles in $\C$ only in poles of $R$, and hence $\Psi = 1/P$ with a polynomial $P$. These poles are of first order and their residues are negative integers since $\Psi = R'/(R-a) - a/(R-a)$. The first term on the right is a logarithmic derivative, while the second term is $0$ at a pole of $R$. Suppose $\Psi$ has two poles in $\C$. Then $\Psi = 1/P$ is $O(1/z^2)$ for $z \to \infty$. This shows that the residue of $\Psi$ at $\infty$ is zero, while the residue theorem implies that it is the negative of the sum of the residues of $\Psi$ in $\C$, which is a positive integer. We get that $\Psi$ has only one pole, so that $\Psi$ has the form $\Psi(z) = -n/(z-p)$ with $n \in \N$. The first order linear differential equation $\Psi(z) = -n/(z-p)$ can be integrated elementarily, which gives the solutions
$$
R(z) = \frac{a}{n+1} (z-p) + a + \frac{C}{(z-p)^{n}}.
$$
$C=0$ is impossible since then $R$ has no pole at $z = p$, so $R$ has to be of the form stated in the theorem.

\glqq$\Leftarrow$\grqq. This can be seen by simple calculations. $\square$

\begin{rem} \rm Since a shared value on the sphere is also a shared value in $\C$, the examples given in Example \ref{rational_IM} are also examples that fit into this section. We were not able to decide if substantially different examples exist. In particular, neither on the sphere nor in the plane we could find examples of the form $R = a(1 + P/(Q'P + P'))$ with $Q' \not \equiv 0$, as suggested by the proof of Lemma \ref{gundersen_representation}.
\end{rem}

\section{$R(w)$ and $\lambda w \cdot R'(w)$ on the sphere} \label{section_R_dR_sphere}

In \cite{Sa1} the first author introduced the notion of {\it shared limit values} for meromorphic functions $f,g : \C \to \Cd$: A value $a \in \Cd$ is a shared limit value for $f$ and $g$ if for all sequences $z_n \to \infty$ it holds 
$$
f(z_n) \to a \Leftrightarrow g(z_n) \to a.
$$
In \cite{Sa2} it was proven that if a meromorphic function $f : \C \to \Cd$ shares three finite limit values with its derivative, then $f$ and $f'$ share all limit values (Theorem 1.3 in \cite{Sa2}). If $f$ is entire it was shown that two finite shared limit values lead to the same conclusion (Theorem 1.4 in \cite{Sa2}).

\begin{ex} \label{ex_entire_slv} \rm The function $f(z) := {\text e}^{2z} + \frac{1}{2}$ is entire and shares with its derivative $f'(z) = 2{\text e}^{2z}$ the limit values $1$ and $\infty$, but no other limit value is shared by $f$ and $f'$. In this sense Theorem 1.4 in \cite{Sa2} is sharp. Note that $f$ is of the form $P({\text e}^{2z})$ with the linear polynomial $P(w) = w + \frac{1}{2}$ and that $f$ and $f'$ also share the values $1$ and $\infty$.
\end{ex}

It is natural to search for examples of the form $f(z) = R({\text e}^{\lambda z})$ with a non-constant rational function $R$, such that $f$ and $f'$ share two finite limit values, in order to show that Theorem 1.3 in \cite{Sa2} is sharp. It seems that all known examples of meromorphic functions sharing two finite values with their derivatives are of this form (see e.g.~\cite{Gu2}, formulas (1) and (2)). The requirements for $R$ in connection with shared limit values are more demanding than for shared values, though: $R(w)$ and $\lambda w \cdot R'(w)$ have to share the two values as mappings from $\Cd$ to $\Cd$ by the Casorati-Weierstra{\ss} theorem, while for shared values it would only be from $\C^*$ to $\Cd$, since ${\text e}^{\lambda z}$ omits $0$ and $\infty$.

\begin{defi} \label{def_dR} \rm Let $R \colon \Cd \to \Cd$ be a non-constant rational function and $\lambda \neq 0$. We define $\partial_{z} R(w) := \lambda w \cdot R'(w)$.
\end{defi}

The following lemma will be used in almost every proof in the rest of the paper. In fact, we need it so many times, that we will not always give explicit reference to it.

\begin{lemma} \label{lemma_dR} Let $R$ be a non-constant rational function.\\

\noindent
i) If $R$ takes the value $c \in \C$ with order $k$ at $w=0$ or $w=\infty$, then $\partial_{z} R$ has a zero of order $k$ at $w=0$ resp.~$w=\infty$, and if $R$ has a pole of order $k$ at $w=0$ or $w=\infty$, then $\partial_{z} R$ also has a pole of order $k$ at $w=0$ resp.~$w=\infty$. In particular, $\partial_{z} R(0)$ and $\partial_{z} R(\infty)$ both lie in $\{ 0, \infty \}$.\\

\noindent
ii) If $R(w)=a \in \C$ with $w \in \C^*$ of order $k$, then $\partial_{z} R(w)=0$ with order $k-1$.\\

\noindent
iii) If $n_{\infty}$ is the number of poles of $R$ in $\C^*$ counted without multiplicities, then $\deg(\partial_{z} R) = d + n_{\infty}$, further $\deg(R - \partial_{z} R) \le d + n_{\infty}$.\\

\noindent
iv) $R \equiv \partial_{z} R$ is equivalent to $R(w) = c \cdot w^n$ with some constant $c \neq 0$, $n \in \Z \setminus \{0\}$ and $\lambda = 1/n$.
\end{lemma}

\pr i) This can easily be verified by expanding $R$ and $\partial_{z} R$ in Laurent series at $w=0$ and $w=\infty$. Statement ii) follows directly from the definition of $\partial_{z} R$. iii) This results immediately by counting the poles of $\partial_{z} R$ and taking i) into account. iv) $R \equiv \partial_{z} R$ is a simple differential equation that can be solved by elementary calculations. $\square$\\

\begin{rem} \rm a) Although we have Lemma \ref{lemma_dR}, ii), $R'$ and $\partial_{z} R$ may have different orders at $w \in \C^*$. Take as an example $R(w) := (w-1)^2$ and $\lambda = 1$. Then $R'(w) = 2(w-1)$ has order one at $w=1/2$, while $\partial_{z} R(w) = 2 w (w-1)$ has order two at $w=1/2$.\\

b) In Lemma \ref{lemma_dR}, iii) it is possible that $\deg(R - \partial_{z} R) < d + n_{\infty}$, e.g. a pole of $R$ at $w=0$ might not be a pole of $R - \partial_{z} R$. A simple example is $R(w) = 1/w$ and $\lambda = -1$, so that $R - \partial_{z} R \equiv 0$.
\end{rem}
\noindent
The next theorem shows that Example \ref{ex_entire_slv} is essentially the only example of its type.

\begin{theo} \label{P_dP_a_sphere}
Let $R \colon \Cd \to \Cd$ be a non-constant rational function. Then the following statements are equivalent.\\

\noindent
i) $R$ has no poles in $\C^*$, and $R$ and $\partial_{z} R$ share a finite value $a$.\\

\noindent
ii) $R$ and $\partial_{z} R$ share a finite value $a$ by CM.\\

\noindent
iii) $R$ has the form
$$
R(w) = a \left( 1 \mp \frac{1}{\lambda d} \right) + c \cdot w^{\pm d}
$$
with $c \neq 0$. 
\end{theo}

\pr ii) $\Rightarrow$ i). Since $a$ is shared by CM, we have $\deg(\partial_{z} R) = d$. Let $n_{\infty}$ be the number of poles of $R$ in $\C^*$, then we always have $\deg(\partial_{z} R) = d + n_{\infty}$, so that $n_{\infty}=0$.\\

i) $\Rightarrow$ iii). We consider
$$
\Psi := \frac{\partial_{z} R - a}{R-a}.
$$
We show that $\Psi$ has no poles. First note that since $R$ has no poles in $\C^*$ Lemma \ref{lemma_dR} shows that $\Psi$ has no poles in poles of $R$.\\

If $a \neq 0$ all $a$-points of $R$ in $\C^*$ are simple, while $R(0)=a$ or $R(\infty)=a$ is impossible by Lemma \ref{lemma_dR}.\\

If $a=0$ then $R$ has no zero in $\C^*$, since otherwise counting the zeros infers $\deg(\partial_{z} P) < \deg(P)$, which is impossible by counting the poles. Hence $R$ can only have zeros or poles at $w=0$ and $w=\infty$ by Lemma \ref{lemma_dR}, so that $\partial_{z} R$ has zeros or poles at $w=0$ and $w=\infty$ of the same order.\\

Hence, in any case, $\Psi$ has no poles at $a$-points of $R$, and therefore $\Psi$ has no poles at all, which shows that $\Psi$ is constant.\\

We may assume that $w = \infty$ is a pole of $R$ of order $m$ (otherwise we use the transformation $w \to 1/w$ and $\lambda \to -\lambda$) and it is then easy to see that $\Psi(\infty) = \lambda m$, so that $\Psi \equiv \lambda m$. This is the following first order differential equation of Euler type
$$
w R'(w) - m R(w) = a \left( \frac{1}{\lambda} - m \right).
$$
The solutions are $R(w) = a(1 - 1/(\lambda m)) + c \cdot w^m$ with $c \in \C$. In particular it follows $m=d$. (If $w=0$ is a pole of order $m$ for $R$ then we have to replace $m$ by $-d$.)\\

iii) $\Rightarrow$ ii) It is easy to check that $R$ and $\partial_{z} R$ of the form given in iii) share the value $a$ by CM. $\square$

\begin{rem} \rm We get Example \ref{ex_entire_slv} for $a=1$, $d=1$, $c=1$ and $\lambda=2$, where $w = {\text e}^{\lambda z}$. Any choice of the parameters yields examples of the form $A + B \cdot {\text e}^{C z}$, so that essentially Example \ref{ex_entire_slv} gives the only nontrivial example of an entire function of the form $f(z) = R({\text e}^{\lambda z})$ that shares a finite non-zero limit value with its derivative. In case of $\lambda d = 1$ it follows $\partial_{z} R = R$, so that all values are shared.
\end{rem}

\begin{ex} \rm There are non-trivial examples for the case that $R$ and $\partial_{z} R$ share on $\Cd$ one finite value $a$ by IM. Define $P(w) = w^2 + w + 3/16 = (w + 1/4)(w + 3/4)$ and set $\lambda=1$, then
$$
R(w) = 1 + \frac{P(w)}{\partial_{z} P(w)} = \frac{48w^2 + 32w + 3}{16 w (2w + 1)},
$$
so that
$$
\partial_{z} R(w) = 1 - \frac{P(w)\,\partial_{z}^2 P(w)}{(\partial_{z} P(w))^2} = \frac{-16w^2 - 12w - 3}{16w(2w + 1)^2}.
$$
$R$ has two simple $1$-points in $w=-1/4$ and $w=-3/4$. These are exactly the $1$-points of $\partial_{z} R$, $w=-3/4$ is a simple $1$-point for $\partial_{z} R$, while $w=-1/4$ has order two.
\end{ex}

The next theorem shows that there are no non-trivial examples for $f(z) = R({\text e}^{\lambda z})$ with a non-constant rational function $R$, such that $f$ and $f'$ share two finite limit values. 

\begin{theo} \label{R_dR_ab_sphere}
Let $R \colon \Cd \to \Cd$ be a non-constant rational function such that $R$ and $\partial_{z} R$ share two finite values, then $R \equiv \partial_{z} R$.  
\end{theo}

\pr First note that by Theorem \ref{P_dP_a_sphere} we may assume that $R$ has a pole in $\C^*$, so that $\deg(\partial_{z} R) > d$.\\

First we assume that the shared values $a$ and $b$ are both non-zero. Neither $a$ nor $b$ can be the value of $R(0)$ or $R(\infty)$, because then $\partial_{z} R(0)=0$ or $\partial_{z} R(\infty)=0$ by Lemma \ref{lemma_dR}, so that $a$ resp.~$b$ would not be shared. Hence all $a$-points and $b$-points $w$ lie in $\C^*$ and are simple by Lemma \ref{lemma_dR}. It follows that the total number of $a$-points and $b$-points is $2 d$. Consider $Q(w) := \partial_{z} R(w)/R(w)$. Since the $a$-points and $b$-points of $R$ are $1$-points of $Q$ this implies $\deg Q \ge 2 d$ or $Q \equiv 1$. Assume $Q \not \equiv 1$. We count the poles of $Q$. We first note that $w=0$ and $w = \infty$ are not poles of $Q$ by Lemma \ref{lemma_dR}. We conclude that all zeros and poles of $R$ lie in $\C^*$ and are simple in order to make $\deg Q \ge 2 d$ possible. This implies that $R(0)$ and $R(\infty)$ are non-zero and finite (and different from $a$ and $b$ as was shown above).\\

We consider
$$
\Psi := \frac{\partial_{z} R - R}{(R-a)(R-b)}.
$$
From the above analysis it is easily seen that $\Psi$ has no poles in $\Cd$. Therefore $\Psi$ is constant.\\

Since $\deg(\partial_{z} R) > d$, at least one of the $a$-points of $\partial_{z} R$ is ramified, because their number is $d$. This means that there exists $w_a \in \C^*$, such that $\partial_{z} R(w_a) = a$ and $(\partial_{z} R)'(w_a)=0$.\\

This and the rule of de L'Hospital show
$$
\Psi(w_a) = \frac{-R'(w_a)}{R'(w_a)(R(w_a)-a) + R'(w_a)(R(w_a)-b)} = \frac{1}{b-a}.
$$
For exactly the same reason there exists $w_b \in \C^*$, such that $\partial_{z} R(w_b) = b$ and $(\partial_{z} R)'(w_b)=0$. This gives
$$
\Psi(w_b) = \frac{-R'(w_b)}{R'(w_b)(R(w_b)-a) + R'(w_b)(R(w_b)-b)} = \frac{1}{a-b}.
$$
But this contradicts that $\Psi$ is constant, and hence $Q \equiv 1$ which means $R \equiv \partial_{z} R$.\\


Now assume that $a=0$ and $b \neq 0$ are the shared values. By Lemma \ref{lemma_dR} $R$ and $\partial_{z} R$ have the same number of zeros and poles in $w=0$ and $w=\infty$ counted with multiplicities. Let $n_0$ be the number of zeros of $R$ in $\C^*$ counted without multiplicities. Since $\partial_{z} R$ is $0$ exactly at the same points with an order that is exactly one less than the order of the zero of $R$, we conclude that $\deg(\partial_{z} R) = d - n_0$. Let now $n_\infty$ be the number of poles of $R$ in $\C^*$ counted without multiplicities. Since $\partial_{z} R$ has its poles exactly at the same points with an order that is exactly one higher than the order of the pole of $R$, we conclude that $\deg(\partial_{z} R) = d + n_\infty$. This is only possible if $n_0 = n_\infty = 0$, i.e.~if $R$ has no zeros and no poles in $\C^*$. Representing $R$ as the quotient of two polynomials it follows $R(w) = c \cdot w^n$ with a constant $c$ and $n \in \Z \setminus \{ 0 \}$. We get $\partial_{z} R(w) = \lambda n c \cdot w^n$. Since $R$ and $\partial_{z} R$ share $b \neq 0$, which has to be taken in $\C^*$, it follows $ \lambda n=1$, so that by Lemma \ref{lemma_dR}, iv) we have $R \equiv \partial_{z} R$. $\square$

\begin{rem} \rm Theorem \ref{R_dR_ab_sphere} shows that the only examples $f(z) = R(w)$ with $w = \exp(\lambda z)$ that share two limit values with $f'$ are $f(z) = c \cdot {\text e}^{z}$. Hence there exists no non-trivial example of this type, in particular the examples in \cite{Gu2} (or Example \ref{example_R_dR}) share two values but not two limit values.
\end{rem} 

\section{$R(w)$ and $\lambda w \cdot R'(w)$ in $\C^*$} \label{section_R_dR_C^*}

If $R$ and $\partial_{z} R$ share a value $a$ in $\C^*$, then $f(z) = R(w)$ with $w={\text e}^{\lambda z}$ shares the value $a$ with its derivative $f'(z)$ for $z \in \C$. It seems to us that all known examples of transcendental meromorphic functions $f \colon \C \to \Cd$ that share values with their derivatives are of this type, in particular the trivial example $f(z) = c \cdot{\text e}^{z}$ which is equivalent to $R \equiv \partial_{z} R$.\\

As already mentioned in the introduction, if $f \colon \C \to \Cd$ shares three finite values with $f'$ then by \cite{Gu} (Theorem 3) and \cite{MuSt1} (Satz 2) this implies $f'=f$, so that the trivial example is the only one in this case. It is still an open question whether this conclusion is still true if only two non-zero values are shared. See \cite{MuSt2}, BEISPIEL B and the comments thereafter.\\

Since $\partial_{z} R$ again plays an important role in this section, we remind of Definition \ref{def_dR} and Lemma \ref{lemma_dR}.

\begin{ex} \label{example_R_dR} \rm It seems that all known examples of meromorphic functions sharing two finite values with their derivatives are up to transformations the examples in \cite{Gu2}, formulas (1) and (2). In both examples one of the shared values is $a=0$. (Otherwise they would provide an answer the question posed above.) We translate the two examples to our notation. The example in formula (2) in \cite{Gu2} can (in principal) be written as
$$
R_{(2)}(w) =\frac{2}{1- C w} 
$$
with $\lambda = -2$ and $C \neq 0$. $R_{(2)}$ and 
$$
\partial_{z} R_{(2)}(w) = \frac{-4Cw}{(1- Cw)^2}
$$ 
share the values $0$ and $1$ in $\C^*$, where $0$ is omitted in $\C^*$ by both functions.\\

The example in formula (1) in \cite{Gu2} can (in principal) be written as
$$
R_{(1)}(w) = \left(\frac{A + BCw}{1+Cw} \right)^2 
$$
with $A = (1-\sqrt{5})/2$, $B = (1+\sqrt{5})/2$, $\lambda = \sqrt{5}/2$ and $C \neq 0$. $R_{(1)}$ and 
$$
\partial_{z} R_{(1)}(w) = \frac{5Cw (A + BCw)}{(1+Cw)^3}
$$ 
share the values $0$ and $1$ in $\C^*$.
\end{ex}
\noindent
We start the section with a theorem that is directly related to entire functions $f$, that share two finite values with their derivative. It will also serve as an argument in later proofs.

\begin{theo} \label{R_dR_ab_no_poles_punctured_plane}
Let $R \colon \C^* \to \Cd$ be a non-constant rational function with no poles in $\C^*$ such that $R$ and $\partial_{z} R$ share two finite values, then $R \equiv \partial_{z} R$.
\end{theo}

\pr The claim follows from Satz 1 in \cite{MuSt1} by considering the entire function $f(z) = R(w)$ with $w = {\text e}^{\lambda z}$. We give an alternative proof for rational functions.\\

Without loss of generality we may assume that $R$ has a pole at $w=\infty$, since otherwise we can apply the transformation $w \rightarrow 1/w$, $\lambda \rightarrow -\lambda$.\\

First we assume that the shared values are both non-zero. We consider $Q := \partial_{z} R - R$. Counting poles and Lemma \ref{lemma_dR} show $\deg(Q) \le d$. We assume $Q \not \equiv 0$. Since every $a$- or $b$-point in $\C^*$ (all of order $1$) leads to a zero of $Q$, and since $w=\infty$ is a pole of $R$, it follows that without loss of generality $R$ has no $a$-points in $\C^*$ since otherwise $\deg(Q) > d$. We conclude that the only $a$-point of $R$ is $w=0$. This gives the form $R(w) = a + c \cdot w^d$, and hence $\partial_{z} R = \lambda d c \cdot w^d$. But then it is easy to check that $\partial_{z} R$ has $a$-points in $\C^*$, so that $a$ cannot be shared by $R$ and $\partial_{z} R$. It follows $Q \equiv 0$, which proves the theorem in this case.\\

Now we assume $a=0$. All zeros of $R$ in $\C^*$ are at least of order $2$. Let $n_0$ be number of zeros of $R$ in $\C^*$, counted without multiplicities. It follows $n_0 \le d/2$. We consider $ L := \partial_{z} R/R$ and assume $L \not \equiv 1$. First note that by Lemma \ref{lemma_dR} $L$ has no poles in $w=0$ or $w=\infty$. Counting the poles of $L$ shows $\deg(L) = n_0$. Let $n_b$ be the number of (simple) $b$-points of $R$ in $\C^*$. Since at these $b$-points $L$ has $1$-points, we conclude $n_b \le n_0 \le d/2$. Hence $R$ has a $b$-point at $w=0$ with an order that is greater or equal than $d/2$, so that by Lemma \ref{lemma_dR} $\partial_{z} R$ has a zero at $w=0$ of the same order. By Lemma \ref{lemma_dR} $\deg(\partial_{z} R) = d$, and all $b$-points of $\partial_{z} R$ lie in $\C^*$. Since these $b$-points have to lie in the at most $d/2$ $b$-points of $R$, the ramification of the $b$-points of $\partial_{z} R$ has to be greater or equal than $d/2$. Summing up the ramification of $\partial_{z} R$ at $w=0$ and the $b$-points in $\C^*$ we get at least $d/2 - 1 + d/2 = d-1$. Since the ramification of $\partial_{z} R$ in its only pole at $w=\infty$ equals $d-1$, we conclude that $\partial_{z} R$ is ramified only in its zero at $w=0$, its $b$-points and its pole at $w=\infty$. Hence all zeros of $\partial_{z} R$ in $\C^*$ are simple, and therefore all zeros of $R$ are of order two, and hence $R=Q^2$ with a polynomial $Q$ that has only simple zeros. This implies $\partial_{z} R = 2 \lambda w Q Q'$. Since the zeros of $Q$ are simple, zeros of $Q'$ are not zeros of $R$. It follows that $Q'$ is constant, and $R$ has the form $R(w) = c(w-A)^2$ with non-zero $c$ and $A$. From $b= R(0)=c A^2$ it follows that $w=2A$ is the only $b$-point of $R$ in $\C^*$. We have $\partial_{z} R(w) = 2 \lambda c w (w-A)$, so that $\partial_{z} R(2A) = 4 \lambda c A^2 = b$ or equivalently $\lambda = 1/4$. Then $(\partial_{z} R(w))' = c(2w-A)/2$ is not zero at $w=2A$ and $\partial_{z} R$ has a second $b$-point that is not a $b$-point of $R$. This contradiction shows $L \equiv 1$. $\square$ \\

\begin{theo} \label{R_dR_aCM_C^*}
Let $R \colon \C^* \to \Cd$ be a non-constant rational function. Then the following statements are equivalent.\\

\noindent
i) $R$ and $\partial_{z} R$ share a finite value $a \neq 0$ by CM.\\

\noindent
ii) $R$ has the form
$$
R(w) = a \left( 1 \mp \frac{1}{\lambda d} \right) + c \cdot w^{\pm d}
$$
with $a \neq 0$ and $c \neq 0$.
\end{theo}

\pr This follows more or less directly from Theorem \ref{P_dP_a_sphere}: Let $n_a$ be the number of $a$-points of $R$ in $\C^*$. Since $\partial_{z} R$ cannot take the value $a$ at $w=0$ or $w=\infty$, the $a$-points of $\partial_{z} R$ lie in $\C^*$ and are simple because $a$ is shared by CM. It follows $\deg(\partial_{z} R) =  n_a \le d$, and hence $\deg(\partial_{z} R)=d$ and $n_a=d$, so that $a$ is shared by CM on $\Cd$. Theorem \ref{P_dP_a_sphere} shows the assertion. $\square$\\

\begin{prop} \label{R_dR_0CM_C^*}
Let $R \colon \C^* \to \Cd$ be a non-constant rational function. Then the following statements are equivalent.\\

\noindent
i) $R$ and $\partial_{z} R$ share the value $a = 0$ by CM.\\

\noindent
ii) $R$ has one of the following three forms\\

a) $R(w)=1/P(w)$ with a polynomial $P$ of degree $d$ with $P(0) \neq 0$, such that $P'$ has zeros only in zeros of $P$ or in $w=0$,\\

b) $R(w)=w^d/P(w)$ with a polynomial $P$ of degree $d$ with $P(0) \neq 0$, such that the polar derivative $d \cdot P(w) - wP'(w)$ of $P$ has zeros only in zeros of $P$ or in $w=0$,\\

c)  $R(w) = c \cdot w^{d}$.
\end{prop}

\pr It is easy to see, that $R$ and $\partial_{z} R$ share the value $a = 0$ by CM, if, and only if, $R$ and $\partial_{z} R$ have no zeros in $\C^*$.\\

If $R(0)=R(\infty)=0$, or $R(0)=0 \wedge R(\infty)=\infty$, or $R(0)=\infty \wedge R(\infty)=0$ then $0$ is shared CM on $\Cd$, so that Theorem \ref{P_dP_a_sphere} with $a=0$ shows that $R$ has the form in ii), c) or a).\\

The remaining two cases are that either $w=\infty$ is the only zero of $R$ and $R(0) \in \C^*$, or $w=0$ is the only zero of $R$ and $R(\infty) \in \C^*$.\\ 

If $w=\infty$ is the only zero it follows $R = 1/P$ with a polynomial $P$ of degree $d$ with $P(0) \neq 0$, and such that $\partial_z R(w) = \lambda w P'(w) / (P(w))^2$ has no zeros in $\C^*$. This is the case if, and only if, $P'$ has zeros only in zeros of $P$ or in $w=0$, and hence ii), a) follows.\\

If $w=0$ is the only zero of $R$ and $R(\infty) \in \C^*$ then we have $R(w)=w^d/P(w)$ with a polynomial $P$ of degree $d$. This gives
$$
\partial_{z} R(w) = \lambda w^d \frac{d P(w) - w P'(w)}{(P(w))^2}.
$$
It follows that $d \cdot P(w) - w P'(w)$ has zeros only in zeros of $P$ or in $w=0$, which proves ii), b).\\

It is easy to show that all three options for $R$ in ii) share $0$ CM with $\partial_{z} R$. $\square$\\

\begin{ex} \rm There are many examples of $R$ and $\partial_{z} R$ that share the value $a = 0$ by CM in $\C^*$. E.g.~with $n, m \in \N$
$$
R(w) = \frac{1}{(w^n-1)^m} \quad \text{and} \quad \partial_{z} R(w) = -\lambda n m \frac{w^n}{(w^n-1)^{m+1}}
$$
or
$$
R(w) = \frac{w^{n \cdot m}}{(w^n-1)^m} \quad \text{and} \quad \partial_{z} R(w) = -\lambda n m \frac{w^{n \cdot m}}{(w^n-1)^{m+1}}.
$$
\end{ex}

\noindent
The next theorem hints in the direction, that if $f \colon \C \to \Cd$ shares with $f'$ a non-zero value by CM and a second value IM, then there should be a uniqueness theorem. At least there exist no counterexamples of the form given in Example \ref{example_R_dR}.

\begin{theo} \label{R_dR_aCMbIM_punctured_plane}
Let $R \colon \C^* \to \Cd$ be a non-constant rational function such that $R$ and $\partial_{z} R$ share two finite values $a \neq 0$ and $b$, where $a$ is shared CM, then $R \equiv \partial_{z} R$.
\end{theo}

\pr Theorem \ref{R_dR_aCM_C^*} shows that $R$ is of the form $R(w) = a \left( 1 \mp 1/(\lambda d) \right) + c \cdot w^{\pm d}$, so that $\partial_{z} R(w) = \pm \lambda d c \cdot w^{\pm d}$.\\

If $b \neq 0$, then it is easy to see that $\partial_{z} R$ has a $b$-point $w_b \in \C^*$, which also has to be a $b$-point of $R$. The equations $\partial_{z} R(w_b)=b$ and $R(w_b)=b$ imply by simple calculation the equation $a \left( 1 \mp 1/(\lambda d) \right) = b \left( 1 \mp 1/(\lambda d) \right)$. Since $a \neq b$ this shows $1 \mp 1/(\lambda d) = 0$, and Lemma \ref{lemma_dR}, iv) shows the claim in this case.\\

Suppose $b=0$. Then it is easy to see that $\partial_{z} R$ has no zero in $\C^*$, and hence $R$ has also no zero in $\C^*$. This is only possible if $a \left( 1 \mp 1/(\lambda d) \right)=0$, and since $a \neq 0$ we get $1 \mp 1/(\lambda d) = 0$, and again Lemma \ref{lemma_dR}, iv) proves the claim. $\square$ \\

Now we consider the case that the value $0$ is shared by CM.

\begin{theo} \label{R_dR_a=0CMbIM_punctured_plane}
Let $R \colon \C^* \to \Cd$ be a non-constant rational function such that $R$ and $\partial_{z} R$ share two finite values $a = 0$ and $b$, such that $a$ is shared CM, i.e.~$0$ is omitted by $R$ and $\partial_{z} R$ in $\C^*$. Then $R \equiv \partial_{z} R$, or 
$$
R(w) = \frac{2b}{1 + C \cdot w^d} \qquad \text{or} \qquad R(w) = \frac{2b \cdot w^d}{w^d + C}.
$$
with $C \neq 0$, and $\lambda = -2/d$ or $\lambda = 2/d$, respectively.
\end{theo}

\pr The claim follows from the result in \cite{Zh} by considering $f(z) = R(w)$ with $w = {\text e}^{\lambda z}$. We give an alternative proof for rational functions.\\

We assume $R \not \equiv \partial_{z} R$. Since $0$ is omitted by $R$ in $\C^*$, it follows that $R(0)=0$ or $R(\infty)=0$, say $R(\infty)=0$. Then $K := R(0) \in \C^*$, since otherwise with Lemma \ref{lemma_dR} $a$ and $b$ are shared on $\Cd$ and Theorem \ref{R_dR_ab_sphere} shows $R \equiv \partial_{z} R$. Hence $\infty$ is the only zero of $R$ of order $d$.\\

Let again $n_{\infty}$ be the number of poles of $R$ in $\C^*$ counted without regard to multiplicity, so that by Lemma \ref{lemma_dR} we have $\deg(\partial_{z} R) = d + n_{\infty}$. Lemma \ref{lemma_dR} shows that $\partial_{z} R$ has a zero of order $d$ in $w = \infty$, and thus $w = 0$ is a zero of $\partial_{z} R$ of order $n_{\infty}$. By Lemma \ref{lemma_dR} the $K$-point of $R$ in $w = 0$ also has order $n_\infty$. \\

First we rule out $K=b$. Assume $K=b$. We count the ramification of $\partial_{z} R$: In $w=0$ the ramification is $n_{\infty}-1$ and in $w = \infty$ it is $d-1$. Since all poles lie in $\C^*$ and there are $n_{\infty}$ of them, we conclude that the ramification of $\partial_{z} R$ in the poles is $d + n_{\infty} - n_{\infty} = d$. In the $d - n_{\infty}$ $b$-points the ramification is $d + n_{\infty} - (d - n_{\infty}) = 2 n_{\infty}$. Summing this up we conclude that the ramification of $\partial_{z} R$ is at least $n_{\infty}-1 + d-1 + d + 2 n_{\infty} = 2d + 3 n_{\infty} - 2$. Since the total ramification of $\partial_{z} R$ is $2 \deg(\partial_{z} R) - 2 = 2d + 2 n_{\infty} - 2$ we get $n_{\infty}=0$. But then $R$ has no poles, and $R$ is constant.\\
 
It follows that $R=K/P$ with $K \neq b$ and a polynomial $P$ with $\deg(P)=d$, $P(0)=1$ and $n_{\infty}$ zeros. Consider 
$$
Q(w) := \frac{\partial_{z} R(w)}{R(w)} = - \lambda \cdot w \frac{P'(w)}{P(w)}.
$$
Since all $b$-points lie in $\C^*$, it follows that $R$ has $d$ simple $b$-points $w_b$ at which $Q(w_b) = 1$, so that $\deg(Q) \ge d$. Counting the poles we get $\deg(Q) = n_{\infty}$, so that $d \le n_{\infty} \le d$, and hence $d = n_{\infty}$. It follows $\deg(\partial_{z} R) = 2d$, and therefore $\partial_{z} R$ has a zero of order $d$ in $w=0$, so that the $K$-point of $R$ at $w=0$ also is of order $d$. It follows that the Taylor expansion of $R$ in $w=0$ starts with the terms
$$
R(w) = K + e_d \cdot w^d + \ldots
$$
and from $K \equiv R \cdot P$ it follows $P(w) = 1 + C \cdot w^d$ with $C \neq 0$ since terms of order $k$ in $P$ with $1 \le k \le d-1$ would lead to a contradiction. Hence we have
$$
R(w) = \frac{K}{1 + C \cdot w^d}.
$$
By elementary calculations using the common $b$-points of $R$ and $\partial_{z} R$ it follows $K=2b$ and $\lambda=-2/d$. The second possible form of $R$ stated in the theorem corresponds to changing $w$ to $1/w$ and $\lambda$ to $-\lambda$. $\square$ \\

In the rest of the paper we give results concerning two shared values of $R$ and $\partial_{z} R$ in $\C^*$ without using the assumption CM. We start with a proposition that is related to the result in \cite{Li}.

\begin{prop} \label{lemma_d_infty_ab} Let $R \colon \C^* \to \Cd$ be a non-constant rational function of degree $d$ such that $R$ and $\partial_{z} R$ share two finite non-zero values $a$ and $b$ and let $n_\infty$ be the number of poles of $R$ in $\C^*$ counted without multiplicities. Then either $n_\infty \ge d/2$ or $R \equiv \partial_{z} R$.
\end{prop}

\pr We assume $R \not \equiv \partial_{z} R$. According to Theorem \ref{R_dR_ab_sphere} we may assume $R(\infty)=a$.\\

Let $m_\infty$ be the order of $R$ at $w=\infty$ and $m_0$ be the order of $R$ at $w=0$. By Lemma \ref{lemma_dR} we always have $\deg(\partial_{z} R) = d + n_\infty$ and $\deg(R - \partial_{z} R) \leq d + n_\infty$. We consider different cases.\\

First we assume $R(0)=b$. The $a$- and $b$-points of $R$ in $\C^*$ are simple and zeros of $R - \partial_{z} R$, and therefore $2d - m_\infty - m_0 \le \deg(R - \partial_{z} R) \le d + n_\infty$ or 
\begin{align} \label{d}
d \leq m_\infty + m_0 + n_\infty.
\end{align}
We consider the ramification of $\partial_{z} R$. It is $d$ in the poles in $\C^*$, $m_\infty -1$ at $w=\infty$, $m_0 - 1$ at $w=0$, $d + n_\infty - (d - m_\infty)$ in the $a$-points in $\C^*$ and $d + n_\infty - (d - m_0)$ in the $b$-points in $\C^*$. Summing up these terms we get $d + 2n_\infty -2 + 2m_\infty +2m_0$. Since the total ramification of $\partial_{z} R$ is $2d + 2n_\infty - 2$ this implies $2m_\infty + 2m_0 \leq d$. With (\ref{d}) this gives $d \le d/2 + n_\infty$, which proves the claim in this case.\\

If $R(0)=a$ the proof is almost identical. Again inequality (\ref{d}) holds. We only have to note that the ramification of $\partial_{z} R$ in this case is $d + n_\infty - (d - m_\infty - m_0)$ in the $a$-points in $\C^*$ and $d + n_\infty$ in the $b$-points in $\C^*$, which leads to exactly the same inequalities.\\

In the case $R(0) \in \C \setminus \{a,b\}$ all $b$-points lie in $\C^*$ and inequality (\ref{d}) improves to
\begin{align*} 
d \leq m_\infty + n_\infty.
\end{align*}
The ramification of $\partial_{z} R$ in the $b$-points in $\C^*$ is now $d + n_\infty - d = n_\infty$. These two changes to the above reasoning exactly cancel out, which gives the same conclusion.\\

The last case is $R(0)=\infty$. Again, since all $b$-points lie in $\C^*$ we get $d \leq m_\infty + n_\infty$. Now the ramification in the poles in $\C^*$ changes to $d - m_0$. Again, these two changes exactly cancel out, which completes the proof. $\square$ \\

Our last result deals with the seemingly very difficult problems connected to $f$ and $f'$ sharing two values by IM. It is only a partial result and our proof is long and divided into nine cases. We give a complete proof for the interested reader in section \ref{prop_d_le_2}.\\

The proposition says that if the shared values are both non-zero, then there is no example in the spirit of Example \ref{example_R_dR} with $d \le 2$, and that if one of the shared values is $0$, then the examples given in Example \ref{example_R_dR} are essentially the only ones, again if $d \le 2$.

\begin{prop} \label{prop_deg2_punctured_plane}
Let $R \colon \C^* \to \Cd$ be a non-constant rational function with degree $d \le 2$, such that $R$ and $\partial_{z} R$ share two finite values $a$ and $b$. Then $R \equiv \partial_{z} R$, or $R$ is (up to simple transformations) one of the two examples given in Example \ref{example_R_dR}. In particular, if $a$ and $b$ are both non-zero, then $R \equiv \partial_{z} R$.
\end{prop}

\section{Remarks and open problems}
\noindent
{\bf Remarks and open problems concerning section \ref{Rational_functions_sphere}:}\\

\nopagebreak
\noindent
If $R$ shares one non-zero value $a$ with $R'$ on the whole sphere, then we were able to prove that $R=a(1 + P/(cP + P'))$ with a polynomial $P$ with only simple zeros and $c \in \C$. With the differential operator $d(P) = cP+P'$ we get $R=a(1 + P/d(P))$ and $R'= a(1 - (P d^2(P)/(d(P))^2)$. Such $R$ and $R'$ share $a$ if and only if every zero of $d^2(P)$ is either a zero of $P$ or of $d(P)$. In all our examples we have $c=0$, i.e.~$R = a \left( 1 + P/P' \right)$.\\ 

{\it Question 1:} Given a polynomial $P$ with only simple zeros. If every zero of $d^2(P)$ is either a zero of $P$ or of $d(P)$, does it follow that $c=0$?\\

{\it Question 2:} Are there any other examples of this type that are substantially different from the ones given in Lemma \ref{gundersen_representation} and Example \ref{rational_IM}?\\

\noindent
{\bf Remarks and open problems concerning section \ref{Rational_functions_plane}:}\\

\nopagebreak
\noindent
If $R$ shares one non-zero value $a$ with $R'$ in the plane, then it is implicit in the proof of Lemma \ref{gundersen_representation} that $R=a(1 + P/(Q P + P'))$ with a polynomial $P$ with only simple zeros and a polynomial $Q$.\\

{\it Question 3:} Are there any examples $R=a(1 + P/(Q P + P'))$ where $R$ and $R'$ share $a$ in the plane with $Q \not \equiv 0$?\\

\noindent
{\bf Remarks and open problems concerning section \ref{section_R_dR_sphere}:}\\

\nopagebreak
\noindent
The problem for two (non-zero) shared values described in the introduction and Theorem \ref{R_dR_ab_sphere} strongly suggest the following question. Probably this problem is even more difficult to solve than the corresponding question for shared values.\\

{\it Question 4:} If a transcendental meromorphic function $f \colon \C \to \Cd$ shares two finite limit values with its derivative, does it follow that $f$ and $f'$ share all limit values?\\

\noindent
{\bf Remarks and open problems concerning section \ref{section_R_dR_C^*}:}\\

\nopagebreak
\noindent
The main question concerning this section is obvious, and can be considered as the most interesting, and also seemingly most difficult, problem that we were not able to solve in this paper.\\

{\it Question 5:} Let $R \colon \C^* \to \Cd$ be a non-constant rational function, such that $R$ and $\partial_{z} R$ share two finite non-zero values. Does it follow that $R \equiv \partial_{z} R$?\\

\section{Proof of Proposition \ref{prop_deg2_punctured_plane}} \label{prop_d_le_2}

We need the following lemma.

\begin{lemma} \label{lemma_ab_R(0)=infty} Let $R \colon \C^* \to \Cd$ be a non-constant rational function such that $R$ and $\partial_{z} R$ share two finite values $a \neq 0$ and $b \neq 0$ with $R \not \equiv \partial_{z} R$. If $R(0)=\infty$ and $R(\infty)=a$, then $a = b \left( 1+ \frac{1}{\lambda d} \right)$ and $R = b \left(1 + P/\partial_{z} P\right)$ with a polynomial $P$ with degree $d$.
\end{lemma}

\pr We can write $R = b + P/Q$ with relatively prime polynomials $P$ and $Q$. Since $w=0$ is a pole of $R$ we conclude that $w$ is a factor of $Q$. The $b$-points of $R$ in $\C^*$ are the zeros of $P$.\\

Since $a-b = P/Q(\infty) \in \C^*$, we conclude that $P$ and $Q$ both have degree $d$ and $R$ has $d$ simple $b$-points in $\C^*$. Direct calculation shows
$$
\partial_{z} R = \lambda \left( \frac{P'}{Q/w} - \frac{wPQ'}{Q^2} \right).
$$ 
In the $d$ zeros of $P$ the function $L := \frac{P'}{Q/w}$ has to take the value $b/\lambda$. From $\deg(L) \le d-1$ it follows $L \equiv b/\lambda$, so that $\partial_{z} P = b Q$, and hence
$$
R = b \left(1 + \frac{P}{\partial_{z} P} \right).
$$ 
Now $P / \partial_{z} P (\infty) = 1/(\lambda d)$ shows the assertion.$\square$\\

\begin{lemma} \label{lemma_01_R(0)=infty} Let $R \colon \C^* \to \Cd$ be a non-constant rational function such that $R$ and $\partial_{z} R$ share $0$ and $1$ with $R \not \equiv \partial_{z} R$, such that $0$ is not omitted by $R$ in $\C^*$. If $R(0)=\infty$, then either $R(\infty)=1$, or $R = 1 + P/\partial_{z} P$ with a polynomial $P$ with degree $d$.
\end{lemma}
\noindent
We omit the proof, as it is almost identical to the proof of Lemma \ref{lemma_ab_R(0)=infty}.\\

\noindent
{\bf Proof of Proposition \ref{prop_deg2_punctured_plane}.} First we treat the case $\deg(R) = 1$. Then $R$ has the form $R(w) = (Aw+B)/(Cw+D)$ with $\Delta := AD-CB \neq 0$. This implies $\partial_{z} R = \Delta \lambda w/(Cw+D)^2$. If one of the shared values is shared CM or is omitted, then the claim follows from Theorem \ref{R_dR_aCMbIM_punctured_plane} and \ref{R_dR_a=0CMbIM_punctured_plane}, respectively. Hence we can assume that both values are taken by $R$ in $\C^*$ once of order one, and by $\partial_{z} R$ at the same points both of order two. It follows that $(\partial_{z} R)'(w) = \Delta \lambda (D-Cw)/(Cw+D)^3$ has two zeros in $\C^*$ which is not the case.\\

Now we assume $\deg(R)=2$ and that $R$ is not one of the two examples from Example \ref{example_R_dR}. The rest of the proof is divided into two cases: {\bf A)}~$a$ and $b$ are both non-zero and {\bf B)}~one of the shared values is $0$, without loss of generality $a=0$ and $b=1$.\\

{\bf A)} First we consider the case that $a$ and $b$ are both non-zero and assume $R \not \equiv \partial_{z} R$. We consider several cases, according to how many poles of $R$ lie in $\C^*$:\\

{\bf A i)} If $R$ has no poles in $\C^*$ then the assertion follows from Theorem \ref{R_dR_ab_no_poles_punctured_plane}.\\

{\bf A ii)} Suppose $R$ has exactly one simple pole in $\C^*$. Then $R$ has a second simple pole in $w=0$ or $w = \infty$, say $R(0)=\infty$. Not all $a$- and $b$-points of $R$ lie in $\C^*$, otherwise Theorem \ref{R_dR_ab_sphere} shows $R \equiv \partial_{z} R$. Hence we can assume $R(\infty)=a$ of order $1$. If the order was $2$ then $a$ would be omitted in $\C^*$ by $R$, and the claim follows from Theorem \ref{R_dR_aCMbIM_punctured_plane}. Hence $R$ has two $b$-points $w=A$ and $w=B$ and a simple pole $w=C$ in $\C^*$. We conclude that $R$ has the form
$$
R(w) = b + (a-b) \frac{(w-A)(w-B)}{w(w-C)}.
$$
We also know from Lemma \ref{lemma_ab_R(0)=infty} that $R = b \left(1 + P/\partial_{z} P\right)$ with a polynomial $P$.\\

Comparing these two representations of $R$ yields $C = (A+B)/2$ and $\lambda = b / (2(a-b))$. With these substitutions it is easy to check that $\partial_{z} R((A+B)/4) = b$, so that $(A+B)/4$ is one of the two $b$-points $A$ and $B$. Without loss of generality we assume $(A+B)/4 = A$ which implies $B=3A$ and $C=2A$. This gives
$$
R(w) = b + (a-b) \frac{(w-A)(w-3A)}{w(w-2A)}
$$
and
$$
\partial_{z} R(w) = A b \frac{w^2 - 3 A w + 3A^2}{w(w - 2A)^2}.
$$
It is easy to check that $w = 3A/2$ is the $a$-point of $R$ in $\C^*$. Since $\partial_{z} R(3A/2) = 2b$ we get $a=2b$. With this substitution it is easy to check that $\partial_{z} R$ has three mutually different $b$-points in $\C^*$, namely $w = 3A/2$, $w = (3/2 + \sqrt{5}/2) A$ and $w = (3/2 - \sqrt{5}/2) A$, and therefore $b$ is not shared by $R$ and $\partial_{z} R$.\\

{\bf A iii)} Suppose $R$ has one pole of order $2$ in $\C^*$. Counting the poles we get $\deg(\partial_{z} R) = 3$ and $\deg(R - \partial_{z} R) = 3$. Since the $a$- and $b$-points in $\C^*$ of $R$ are zeros of $R - \partial_{z} R$ it follows that not all $a$- and $b$-points lie in $\C^*$, say $R(\infty)=a$. Since all $a$-points of $\partial_{z} R$ are in $\C^*$ by Lemma \ref{lemma_dR}, $R$ has one $a$-point in $\C^*$. This has to be an $a$-point of $\partial_{z} R$ of order $3$. Counting the ramification of $\partial_{z} R$ in the pole and the $a$- and $b$-points we get $5$. But the total ramification of $\partial_{z} R$ is $4$, a contradiction.\\

{\bf A iv)} Suppose $R$ has two simple poles in $\C^*$. Counting the poles we see $\deg(\partial_{z} R) = 4$, so that the total ramification of $\partial_{z} R$ is $6$. If $R$ has only one $a$-point in $\C^*$, then this is an $a$-point of $\partial_{z} R$ of order $4$. Counting the ramification of $\partial_{z} R$ in the poles and in the $a$- and $b$-points gives at least $7$, which is impossible. Likewise, $R$ cannot have only one $b$-point in $\C^*$. Hence all $a$- and $b$-points of $R$ lie in $\C^*$. But then $a$ and $b$ are shared with respect to $\Cd$ and the claim follows from Theorem \ref{R_dR_ab_sphere}.\\

{\bf B)} Now we assume $a=0$ and $b=1$. If $R$ has no zero in $\C^*$, then the claim follows from Theorem \ref{R_dR_a=0CMbIM_punctured_plane}, hence we may assume that $R$ has a zero in $\C^*$.\\

We make case distinctions according to how many poles of $R$ lie in $\C^*$.\\

{\bf B i)} If $R$ has no poles in $\C^*$ then the assertion follows from Theorem \ref{R_dR_ab_no_poles_punctured_plane}.\\

{\bf B ii)} If $R$ has exactly one simple pole in $\C^*$ then $R$ has a second simple pole in $w=0$ or $w = \infty$, without loss of generality in $w=0$. If $R(\infty)=0$, then by counting the zeros we get $\deg(\partial_{z} R) < d$, which is impossible. It follows that $R$ has exactly one double zero in $\C^*$. Lemma \ref{lemma_01_R(0)=infty} shows that either $R(\infty)=1$, or $R = 1 + P/\partial_{z} P$ with a polynomial $P$ with degree $2$.\\

Counting the poles and using Lemma \ref{lemma_dR} we get $\deg(\partial_{z} R)=3$. Since $\partial_{z} R$ has only one simple zero in $\C^*$ the zero of $\partial_{z} R$ at $w=\infty$ has order $2$, so that by Lemma \ref{lemma_dR} the order of $R$ at $w=\infty$ is also $2$. If $R(\infty)=1$ then this is the only $1$-point of $R$. But all $1$-points of $\partial_{z} R$ lie in $\C^*$, so that $1$ cannot be a shared value.\\

It follows $R = 1 + P/\partial_{z} P$ with a polynomial $P(w) = A + Bw + w^2$. Direct calculation shows that $R$ has order two at $w = \infty$ if and only if $B=0$. But then it is easy to see that $R$ has a double pole in $w=0$ and no pole in $\C^*$, a contradiction.\\

{\bf B iii)} Suppose $R$ has one pole of order $2$ in $\C^*$. Since $R$ is not the example $R_{(2)}$ from Example \ref{example_R_dR} we deduce from Theorem \ref{R_dR_a=0CMbIM_punctured_plane} that $R$ has a zero in $\C^*$ of order $2$. Hence $R$ has the form $R(w) = (M(w))^2$ with a M{\"o}bius transformation $M(w) = (Aw + B)/(Cw + D)$ where $A \neq 0$ and $C \neq 0$. Without loss of generality we assume $C=1$ and we have
$$
\partial_{z} R(w) = 2 \lambda w \frac{(Aw + B)(AD - B)}{(w + D)^3}.
$$

Clearly $R$ has $1$-points in $\C^*$ exactly at the $1$- and $(-1)$-points of $M$.\\

First we exclude the possibility that $M$ is $1$ or $-1$ in $w=0$ or $w = \infty$. Without loss of generality we assume $M(\infty)=-1$, which is equivalent to $A=-1$. The other cases can be obtained by changing $M$ to $-M$ and/or $\lambda$ to $-\lambda$. The equation $1 = M(u)$ then has the unique solution $u = (B-D)/2$, so that in particular $B \neq D$, since otherwise $R$ omits $1$ in $\C^*$. It follows that $u$ is the only $1$-point of $\partial_{z} R$, which must be of order three since $\deg(\partial_{z} R)=3$, so that $(\partial_{z} R)'(u)=0$ and $(\partial_{z} R)''(u)=0$. The requirement $(\partial_{z} R)'(u)=0$ implies $B=5D/3$. But with this substitution $(\partial_{z} R)''(u)= (63 \lambda)/(8D^2) \neq 0$, a contradiction.\\

Hence we can assume that the $1$- and $-1$-point of $M$ both lie in $\C^*$, in particular we have $A \neq \pm 1$. The equation $M(w)=1$ has the unique solution $u=(D-B)/(A-1)$ and $M(w)=-1$ the solution $v=-(D+B)/(A+1)$. From $1 = \partial_{z} R(u)$ we get by direct calculation $\lambda = (DA-B)/(2(D-B)(A-1))$ and from $-1 = \partial_{z} R(v)$ in a similar fashion $\lambda = (DA-B)/(2(D+B)(A+1))$. Since $\lambda$ has a fixed value, we get by comparing these two expressions $D=-AB$, and therefore $\lambda = (A^2 + 1)/(2(A^2 - 1))$. One of the two $1$-points has to be double for $\partial_{z} R$, say $(\partial_{z} R)'(u)=0$. Solving this equation by direct calculation shows $A = (1 \pm \sqrt{5})/2$. The parameter $B$ can be chosen to be any non-zero value, leading to
$$
R(w) = \left( \frac{\frac{1 \pm \sqrt{5}}{2} w + B}{w - \frac{1 \pm \sqrt{5}}{2} B} \right)^2.
$$
But this is basically the example $R_{(1)}$ from Example \ref{example_R_dR}.\\

{\bf B iv)} We assume that $R$ has two simple poles in $\C^*$, so that $R$ has the form
$$
R(w) = c \cdot \frac{(w-A)^2}{(w-B)(w-C)}
$$
and
$$
\partial_{z} R(w) = \lambda c \cdot \frac{w (w-A) ((2A-B-C)w +2BC - AB - AC)}{(w-B)^2 (w-C)^2}.
$$
Since $0$ is shared, $\partial_{z} R$ has only one simple zero in $\C^*$ at $w=A$. This is only possible if $2A-B-C=0$, which is equivalent to $A=(B+C)/2$, or $2BC - AB - AC=0$, which is equivalent to $1/A=(1/B + 1/C)/2$. The two cases correspond to each other when the above reasoning is applied to
$$
\widetilde{R}(w) = R\left( \frac{1}{w} \right) = c \cdot \frac{\left(w-\frac{1}{A}\right)^2}{\left(w-\frac{1}{B}\right)\left(w-\frac{1}{C}\right)}.
$$ 
Thus we can assume without loss of generality that $A=(B+C)/2$. Since $A \neq 0$ we get in particular $B + C \neq 0$. This gives
$$
R(w) = c \cdot \frac{(w-(B+C)/2)^2}{(w-B)(w-C)}
$$
and
$$
\partial_{z} R(w) = -\lambda c \frac{(B-C)^2}{2}\cdot \frac{w (w-(B+C)/2)}{(w-B)^2 (w-C)^2}.
$$
Suppose $R(\infty)=1$. Then $c=1$, and since $\partial_{z} R$ has a zero of order two at $w = \infty$, $R$ also has order two at $w = \infty$, so this is the only $1$-point of $R$, and hence $R$ omits $1$ in $\C^*$ and the claim follows from Theorem \ref{R_dR_aCMbIM_punctured_plane}.\\
 
We prove $R(0) \neq 1$. Again we consider $\widetilde{R}(w) = R( 1/w )$ which shares $0$ and $1$ in $\C^*$ with $\partial_{z} \widetilde{R}$ with $\widetilde{\lambda} = -\lambda$. The same reasoning as above shows $\widetilde{R}(\infty) \neq 1$, which is equivalent to $R(0) \neq 1$.\\

Hence, $R$ has two simple $1$-points in $\C^*$, say in $w=u$ and $w=v$. The equation $R(w)=1$ leads to the quadratic equation
$$
w^2 - (B+C)w + K = 0
$$ 
with a constant $K$ that will be of no importance. It follows $u+v=B+C$.\\

The equation $\partial_{z} R(w) = 1$ leads to the quartic equation
$$
w^4 - 2(B+C)w^3 + \text{lower terms} = 0.
$$

There are two cases to consider: $\partial_{z} R$ has order three in $u$ and order one in $v$, or the order of $\partial_{z} R$ is two in both $u$ and $v$.\\

If $R$ has order three in $u$ and order one in $v$, then we get $3u + v = 2(B+C) = 2u + 2v$, which gives the contradiction $u=v$. If the order of $R$ at $u$ and $v$ is both two, then consider $Q := R(2 - R) - \partial_{z} R$. Counting poles gives $\deg(Q) \le 4$, but $Q$ has five zeros, when counted according to multiplicity, namely (at least) double zeros at $u$ and $v$ and a further zero at the common zero $w=(B+C)/2$ of $R$ and $\partial_{z} R$. Hence $R(2 - R) \equiv \partial_{z} R$. But this is impossible, since then $\partial_{z} R$ also has a zero of order two in $w=(B+C)/2$, which is not the case.\\

With this, the proof of B~iv), and therefore of Proposition \ref{prop_deg2_punctured_plane} is complete. $\square$

\end{document}